\documentclass{article}
\usepackage{latexsym}
\usepackage{verbatim}
\usepackage{amsmath}
\usepackage{amssymb}
\usepackage{mathrsfs}
\usepackage{amsthm}
\def\Ind#1#2{#1\setbox0=\hbox{$#1x$}\kern\wd0\hbox to 0pt{\hss$#1\mid$\hss}
\lower.9\ht0\hbox to 0pt{\hss$#1\smile$\hss}\kern\wd0}
\def\dnf{\mathop{\mathpalette\Ind{}}}
\def\Notind#1#2{#1\setbox0=\hbox{$#1x$}\kern\wd0\hbox to 0pt{\mathchardef
\nn=12854\hss$#1\nn$\kern1.4\wd0\hss}\hbox to
0pt{\hss$#1\mid$\hss}\lower.9\ht0 \hbox to
0pt{\hss$#1\smile$\hss}\kern\wd0}

\def\ind#1#2#3{{  { #1 {{\dnf}_{#3}} #2 } }}
\begin{document}

\def\Aut{\mathop{\rm Aut}\nolimits}
\def\Sym{\mathop{\rm Sym}\nolimits}
\def\Max{\mathop{\rm Max}\nolimits}
\def\an{\mathop{\rm an}\nolimits}
\def\SL{\mathop{\rm SL}\nolimits}
\def\eq{\mathop{\rm eq}\nolimits}
\def\Ker{\mathop{\rm Ker}\nolimits}
\def\Fix{\mathop{\rm Fix}\nolimits}

\def\Cl{\mathop{\rm Cl}\nolimits}
\def\val{\mathop{\rm val}\nolimits}
\def\Capa{\mathop{\rm Cap}\nolimits}
\def\Min{\mathop{\rm Min}\nolimits}
\def\prof{\mathop{\rm prof}\nolimits}
\def\Th{\mathop{\rm Th}\nolimits}
\def\lim{\mathop{\rm lim}\nolimits}
\def\Sup{\mathop{\rm Sup}\nolimits}
\def\rk{\mathop{\rm rk}\nolimits}
\def\gp{\mathop{\rm gp}\nolimits}
\def\Fitt{\mathop{\rm Fitt}\nolimits}
\def\Soc{\mathop{\rm Soc}\nolimits}
\def\PSL{\mathop{\rm PSL}\nolimits}
\def\SU{\mathop{\rm SU}\nolimits}
\def\rings{\mathop{\rm rings}\nolimits}
\def\Alt{\mathop{\rm Alt}\nolimits}
\def\prof{\mathop{\rm prof}\nolimits}
\def\an{\mathop{\rm an}\nolimits}
\def\tp{\mathop{\rm tp}\nolimits}

\newtheorem{defi}{Definition}[subsection]
\newtheorem{theorem}[defi]{Theorem}
\newtheorem{definition}[defi]{Definition}
\newtheorem{lemma}[defi]{Lemma}
\newtheorem{proposition}[defi]{Proposition}
\newtheorem{conjecture}[defi]{Conjecture}
\newtheorem{corollary}[defi]{Corollary}
\newtheorem{problem}[defi]{Problem}
\newtheorem{remark}[defi]{Remark}
\newtheorem{example}[defi]{Example}
\newtheorem{question}[defi]{Question}
\newtheorem{convention}[defi]{Conventions}
\newtheorem{fact}[defi]{Fact}

\title{Model theory of finite and pseudofinite  groups}
\author{Dugald Macpherson\footnote{Research partially supported by EPSRC grant EP/K020692/1},\\ School of Mathematics,\\University of Leeds,\\ Leeds LS2 9JT, UK,\\h.d.macpherson@leeds.ac.uk}

\maketitle


\begin{abstract}
 This is a survey, intended both for group theorists and model theorists,  concerning the structure of {\em pseudofinite groups}, that is, infinite models of the first order theory of finite groups. The focus is on concepts from stability theory and generalisations in the context of pseudofinite groups, and on the information this might provide for finite group theory.
\end{abstract}

\section{Introduction}
This article is mainly a survey, based on notes for a lecture course at the `Models and Groups 5' meeting in Istanbul October 8--10 2015, but closely related to material on pseudofinite structures which I discussed in the `IPM conference on set theory and model theory', Tehran, October 12-16 2015. The focus below is mainly on pseudofinite groups which are simple in the group-theoretic sense, on the content  for pseudofinite groups of model-theoretic tameness conditions generalising stability, and on the implications for finite group theory. The paper is intended for both logicians and group theorists, so contains considerably more  model-theoretic background than is standard for an article in a logic journal.

\medskip

{\bf Convention}: We let $L_{\gp}:=(\cdot, {}^{-1},1)$ be the first order language of groups. Unless otherwise mentioned, any first order language $L$ is assumed to be countable.

\begin{definition} A {\em pseudofinite group} is an infinite group which satisfies every first order sentence of $L_{\gp}$  that is true of all finite groups.
\end{definition}

Not every group is pseudofinite. For example, the sentence (for abelian groups, so written additively) expressing
`if the map $x\mapsto 2x$ is injective then it is surjective' is true in all finite groups but false in $({\mathbb Z},+)$, and so the latter is not pseudofinite. Likewise (considering the map $x \mapsto px$),
the group of $p$-adic integers $({\mathbb Z}_p,+)$ is not pseudofinite. Since centralisers of non-identity elements in free groups are definable and isomorphic to $({\mathbb Z},+)$, free groups are not pseudofinite. Answering a question asked in Istanbul by G. Levitt, we show in Theorem~\ref{levittq} that there is a group (namely the full symmetric group on a countably infinite set) which does not embed in any pseudofinite group. 

\begin{remark} \rm A group $G$ is pseudofinite if and only if it is elementarily equivalent to a non-principal ultraproduct (see Section 2) of distinct finite groups.

In fact, the above definition, and this remark, make sense with `group' replaced by `field', `ring', `graph', $L$-structure, etc. In this paper, we also talk of pseudofinite rings and fields.
\end{remark}

The structure of this article is as follows. Section 2 contains an overview of basic background, around ultraproducts, pseudofinite fields, and basic concepts from generalised stability theory. In Section 3, we discuss three major theorems of John Wilson about finite and pseudofinite groups: the description of simple pseudofinite  groups; the finite axiomatisability, among finite groups, of soluble groups; and the uniform definability of the soluble radical of a finite group. In Section 4 we consider which pseudofinite groups have a first order theory which is stable, or simple or NIP, or NTP${}_2$  -- this last appears to be new, though straightforward. In Section 4 we also discuss the notions of {\em asymptotic class} of finite structures,  introduced by Elwes, myself and Steinhorn, and consider  these  in the context of groups. In Section 5 we discuss  pseudofinite permutation groups, especially material from \cite{lmt}. We take the opportunity here  to correct some inaccuracies in \cite{ejmr} and \cite{lmt}. Section 6 considers some further directions, and contains a small new result.

\medskip

{\em Acknowledgement:} I warmly thank the organisers of the `IPM conference on set theory and model theory', Tehran, October 12-16 2015, a meeting which led to preparation of this paper. This work was partially supported by EPSRC grant EP/K020692/1.

\section{Preliminaries}

\subsection{Ultraproducts.} 
 Fix a countable language $L$ -- the undelying languages considered throughout this paper are countable. Let $\{M_i:i\in \omega\}$ be a family of $L$-structures, and let $\mathcal{U}$ be a non-principal ultrafilter on $\omega$. (An {\em ultrafilter} on $\omega$ is a family of subsets of $\omega$ closed under finite intersections and supersets, containing $\omega$ and omitting $\emptyset$, and maximal subject to this; it is {\em principal} if it has the form $\{X\subseteq \omega:a\in X\}$ for some $a\in \omega$, and is {\em non-principal} otherwise.) Define $M^*:=\Pi M_i$ (the Cartesian product of the $M_i$.) We say that some property $P$ holds {\em almost everywhere} or {\em for almost all $i$} if 
$$\{i: P \mbox{~holds for~} M_i\} \in \mathcal{U}.$$ For $a=(a_i)_{i\in \omega}$ and
$b=(b_i)_{i\in \omega}$, put $a\sim b$ if $\{i:a_i=b_i\}\in \mathcal{U}$. Then $\sim$ is an equivalence relation. Put $M=M^*/\sim$. Define relations of $L$ to hold of a tuple of $M$ if they hold in the $i$ coordinate (that is, in $M_i$) for almost all $i$, and interpret functions and constants in $M$ similarly.  This is well-defined, and the resulting $M$ is called the ultraproduct of the $M_i$ with respect to $\mathcal{U}$, and here denoted $\Pi_{i\in \omega} M_i/\mathcal{U}$. 
The ultraproduct $M$ will be $\omega_1$-saturated: any type over any countable subset of $M$ will be realised in $M$.
The key fact about ultraproducts is 

\begin{theorem}[{\L}os's Theorem]
In the above notation, for any sentence $\sigma$, $M\models \sigma$ if and only if $\sigma$ holds of $M_i$ for almost all $i$.
\end{theorem}

The following well-known observation makes a link to pseudofiniteness.

\begin{proposition}
A group (or just $L$-structure) is pseudofinite if and only if it is elementarily equivalent to   an infinite ultraproduct of finite structures.
\end{proposition}

\subsection{Pseudofinite fields.} 
We summarise aspects of the  beautiful  structure theory of pseudofinite {\em fields}. This originated with Ax in 1968, and is essential for understanding pseudofinite groups which are simple as groups (see Theorem~\ref{simplegroup} below).

\begin{theorem}[\cite{ax}]
A field $F$ is pseudofinite if and only if all of the following hold: 

(i) $F$ is perfect;

(ii) $F$ is {\em quasifinite} (that is, inside a fixed algebraic closure, $F$ has a unique extension of each finite degree); 

(iii) $F$ is {\em pseudo-algebraically closed} (PAC), that is, every absolutely irreducible variety which is defined over $F$ has an $F$-rational point.

\end{theorem}

It is easily seen that (i) and (ii) hold of all finite fields, and are first-order expressible ((ii) needs some work). (iii) is expressible by a conjunction of first order sentences (this is not completely obvious) each of which, by the Lang-Weil estimates, holds in sufficiently large finite fields, and so each must hold of any pseudofinite field. The striking fact is the converse, that any field satisfying all of these three conditions satisfies {\em every} sentence true of all finite fields. 

Ax also identified the complete theories of pseudofinite fields. If $F$ is a field, then Abs$(F)$ denotes the intersection of $F$ with the algebraic closure of its prime subfield.

\begin{theorem}[\cite{ax}] If $F_1,F_2$ are pseudofinite fields, then $F_1\equiv F_2$ (that is, they are elementarily equivalent) if and only if $F_1,F_2$ have the same characteristic and {\rm Abs}$(F_1)\cong${\rm Abs}$(F_2)$.
\end{theorem}

This, with further information in \cite{ax}, was used by Kiefe \cite{kiefe} to prove a uniform partial quantifier elimination in finite fields, and hence for the theory of  pseudofinite fields: any formula $\phi(\bar{x})$ in the language $L_{\rings}$ 
of rings is equivalent, modulo the theory of finite fields, to a boolean combination of sentences of the form $\exists y g(\bar{x},y)=0$, where $g(\bar{X},Y)\in {\mathbb Z}[\bar{X},Y]$. This can be converted into a model completeness result after the language is
 expanded by constants (see \cite{cdm}). It is also known (see \cite[corollary 3.1]{hrushrav}) that any complete theory of pseudofinite fields has elimination of imaginaries over constants naming an elementary submodel.
 
 \subsection{Basics of stability theory and generalisations.}

 We will consider in this paper {\em stable} theories of pseudofinite groups,  the orthogonal generalisations {\em simple} and {\em NIP} of stable, and their common generalisation {\em NTP${}_2$}. Here, we briefly introduce these model-theoretic concepts. There are many sources on stability theory -- see for example \cite{pillay}, or \cite{wagner} for stable groups. For background on simple theories see \cite{wagner2}, \cite{casan}, or \cite{kim}, and for NIP theories see \cite{simon}. An excellent source of general background is \cite{tz}.

Below, given a complete theory $T$, we let $\bar{M}$ denote a `sufficiently saturated' model of $T$, with all parameter sets taken inside $\bar{M}$.

\begin{definition} \rm 
Let $T$ be a complete theory. A formula $\phi(\bar{x},\bar{y})$ is {\em unstable} (for $T$) if there are $\bar{a}_i\in \bar{M}^{|x|}$ and $\bar{b}_i\in \bar{M}^{|y|}$ (for all $i\in \omega$) such that for any $i,j\in \omega$,
$\bar{M}\models \phi(\bar{a}_i,\bar{b}_j)$ if and only if $i<j$.

The theory $T$ is {\em stable} if no formula is unstable for $T$.
\end{definition}

Several other conditions are equivalent to stability. For example, for  $A\subset \bar{M}$ let
$S_n(A)$ be the set of all $n$-types over $A$. Then $T$ is {\em $\lambda$-stable} (for $\lambda$ an infinite cardinal) if for all $A\subset \bar{M}$ with $|A|\leq \lambda$ we have $|S_1(A)|\leq \lambda$, and $T$ is stable if and only if it is $\lambda$-stable for some infinite $\lambda$.

A theory $T$ is stable if and only if there is an `independence relation' $\ind{A}{B}{C}$  (read `$A$ is independent from $B$ over $C$') satisfying a number of natural axioms 
(suggested by linear independence in vector spaces, or algebraic independence in fields) such as  {\em symmetry}: $\ind{A}{B}{C} \Leftrightarrow \ind{B}{A}{C}$.
One of these axioms is 

{\em local character}: for any $\bar{a}$ and $B$ there is countable $B_0\subset B$ such that
$\ind{\bar{a}}{B}{B_0}$ (here we assume the underlying language is countable).

 Another is {\em stationarity}: 

if $A\subset B\subset M$ and $A$ is algebraically closed in the sense of $T^{\eq}$, and $\bar{a}\in \bar{M}$, then there is $\bar{a}'\in \bar{M}$ such that $\tp(\bar{a}/A)=\tp(\bar{a}'/A)$, and $\ind{\bar{a}'}{B}{A}$. 
 
 \noindent In a stable theory, the independence is given by {\em non-forking} (not defined here).
 
 \begin{definition} \rm
 A formula $\phi(\bar{x},\bar{y})$ has the {\em tree property} (with respect to $T$) if for some $k\in \omega$ the following hold: there are $\bar{a}_\eta\in \bar{M}^{|\bar{y}|}$ 
for all $\eta \in {}^{<\omega}\omega$ such that  for any $\eta \in {}^{<\omega}\omega$ the
set
 $\{\phi(\bar{x},\bar{a}_{\eta i}):i\in \omega\}$ is $k$-inconsistent (that is, any intersection of size $k$ is inconsistent), and for any $\sigma \in {}^\omega\omega$,
 the set
 $\{\phi(\bar{x},\bar{a}_\eta): \eta \mbox{~restricts~} \sigma\}$ is consistent. 
 
 The theory $T$ is {\em simple} if no formula has the tree property.
 \end{definition}
 
 There is a characterisation of simplicity like the above one for stability, via an independence relation $\downarrow$, with the `stationarity' axiom  weakened to the `independence theorem',
 also called `type amalgamation'. Simplicity is a proper generalisation of stability. Within the class of simple theories is that of {\em supersimple} theories, characterised among simple theories by
 a strengthening of the local character condition on $\downarrow$: a simple theory $T$ is {\em supersimple} if and only if, given any $\bar{a}$ and $B$, there is {\em finite} (as distinct from just countable) $B_0\subseteq B$ such 
that $\bar{a}\downarrow _{B_0} B$. For supersimple theories, there is a notion of ordinal-valued {\em rank} on definable sets (or types), known as SU-rank, which we do not here define.
 
 \begin{definition} \rm  A formula $\phi(\bar{x},\bar{y})$ has the {\em independence property} (for $T$)
 if there are $M\models T$ and $\bar{a}_i\in M^{|\bar{x}|}$ for each $i\in \omega$ 
such that for any $S\subset \omega$ there is $\bar{b}_S\in M^{|\bar{y}|}$ with, for each $i\in \omega$,  
$M\models \phi(\bar{a}_i,\bar{b}_S)$ if and only if $i\in S$.

A complete theory $T$ has the {\em independence property} if some formula has the independence property for $T$. We say $T$ is {\em NIP} if it  does not have the independence property. NIP theories are also called {\em dependent} theories. 
\end{definition}

\begin{example} \label{exstable} \rm 
Examples of $\omega$-stable theories include algebraically closed fields, and (hence), for an algebraically closed field $K$, the $K$-rational points of an algebraic group defined over $K$. Separably closed fields which are not algebraically closed are stable but not $\omega$-stable. Abelian groups (and more generally, modules, in the usual language of modules over a fixed ring) are stable, as are free groups.

Any o-minimal structure is  NIP but not stable, as is ${\mathbb Q}_p$,  any non-trivially valued algebraically closed field (in a language defining the valuation), and many other henselian valued fields.

Pseudofinite fields are not stable. For example, if $F$ is a pseudofinite field of odd characteristic, and $\phi(x,y)$ is the formula $\exists z(z^2=x-y)$, then $\phi$ has the independence property. However, pseudofinite fields have simple theory. In fact, they are supersimple  of  SU-rank 1. The well-known theory ACFA (the model companion of the theory of  fields equipped with an automorphism) has all its completions  supersimple, of SU-rank $\omega$; such a field is algebraically closed, and the fixed field of the automorphism is pseudofinite. Groups such as PSL${}_n(F)$ (where $F$ is a pseudofinite field) will have supersimple finite rank theory, but are unstable because they interpret the underlying field $F$.
\end{example}

Suppose that $G$ is a group definable in an NIP theory $T$, and let $\phi(x,\bar{y})$ be any formula. By the Baldwin-Saxl Theorem (\cite{bs}, see also \cite{wagner}), there is $n_\phi\in \omega$ such that any finite intersection of $n_\phi$ $\phi$-definable subgroups of $G$
(i.e. a subgroup of form $\bigcap_{i=1}^t \phi(G,\bar{a}_i)$, where the $\phi(G,\bar{a}_i)$ are subgroups of $G$) is an intersection of at most $n_\phi$ of them. If in addition $T$ is {\em stable}, then (essentially because $T$ cannot have the `strict order property'), this ensures that $G$ has the descending chain condition on intersections of $\phi$-definable subgroups of $G$ -- there is a fixed bound on the lengths of such chains. In particular, we may apply this to the formula $\phi(x,y)$ expressing $xy=yx$. If $T$ is NIP then there is $n_\phi$ such that for any {\em finite} $F\subset G$ there is $F_0\subset F$ with $|F_0|\leq n_{\phi}$ such that $C_G(F)=C_{G}(F_0)$, and if in addition $G$ is stable then {\em any} chain of centralisers has bounded length.

Finally, a complete theory $T$ is {\em TP${}_2$} (has the {\em tree property of the 2nd kind}) if there are $\{\bar{b}_{ij}:i,j< \omega\}$ in $\bar{M}\models T$ and $k<\omega$ such that

(i) the set $\{\phi(\bar{x},\bar{b}_{ij}): j<\omega\}$ is $k$-inconsistent for each $i<\omega$, and 

(ii) for all $\xi\in \omega^\omega$, the set $\{\phi(\bar{x},\bar{b}_{i,\xi(i)}): i\in \omega\}$ is consistent. 

The theory  $T$ is {\em NTP${}_2$} is it is not TP{${}_2$.

It is known (see \cite{chern1} that in the above definition, we may take $|\bar{x}|=1$ . Any simple or NIP theory is NTP${}_2$. Examples of structures whose theory is NTP${}_2$ but not simple or NIP include: non-principal ultraproducts (over $p$) of fields ${\mathbb Q}_p$, and the universal homogeneous ordered graph. For groups, we have the following useful result of Chernikov, Kaplan, and Simon.

\begin{proposition} \cite{cks} \label{ntp2crit}
Let $T$ be NTP${}_2$, let $G$ be a definable group in $M\models T$, and let $(H_i)_{i\in \omega}$ be  uniformly definable normal subgroups of $G$. Let $H:=\bigcap_{i\in \omega} H_i$, and put $H_{\neq j}:=\bigcap_{i\in \omega\setminus \{j\}} H_i$. Then there is some $i^*\in \omega$ such that $|H_{\neq i^*}:H|$ is finite.
\end{proposition}

 \section{Three theorems of  Wilson}

We  consider first {\em simple} groups which are pseudofinite. We warn the reader that in this paper we consider both simple groups (groups with no proper non-trivial normal subgroups) and simple theories (complete theories  for which no formula has the tree property), and that the word `simple' may have both meanings in the same sentence. For background on groups of Lie type, including twisted groups, see for example Carter \cite{carter}. The groups of Lie type are determined by a Dynkin diagram, and a field, and (for the twisted groups) a symmetry of the Dynkin diagram. 

\begin{theorem}[Wilson \cite{wilson1}] \label{simplegroup}
A pseudofinite group $G$ is simple if and only if $G$ is a simple group of Lie type (possibly twisted) over a pseudofinite field.
\end{theorem}

\begin{remark} \rm
1.  In \cite{wilson1} the statement is just that $G$ is {\em elementarily equivalent} to such a group of Lie type; the assertion as given uses also work of Ryten \cite{ryten} discussed in Section 4.1.

2. Ugurlu \cite{ugurlu} has shown that one can replace `simple' by `definably simple of finite centraliser dimension'. Here, a group is {\em definably simple} if it has no proper non-trivial {\em definable} normal subgroups. We say that $G$ has {\em centraliser dimension $k$} if $k$ is the largest natural number such that there is a sequence
$$G=C_G(x_0)>C_G(x_0,x_1)>\ldots > C_G(x_0,\ldots,x_k)=Z(G),$$
and $G$ has {\em finite centraliser dimension} if $G$ has centraliser dimension $k$ for some natural number $k$.

\end{remark}

The proof of the direction $\Leftarrow$ of Theorem~\ref{simplegroup} follows from the fact that finite simple groups of fixed Lie type $\tau$ are {\em boundedly simple}: there is $d=d(\tau)\in \omega$ such that if $G$ is such a group and $g,h\in G$ with $h\neq 1$, then $g$ is a product of at most $d$ conjugates of $h$ and $h^{-1}$. It uses the following result of Point \cite{point}.

\begin{theorem}[Point]
Let $\{G(q_i): i\in I\}$ be a family of finite simple groups of the same Lie type (possibly twisted), and let $\mathcal{U}$ be a non-principal ultrafilter on $\omega$. Then
$$\Pi|_{i\in \omega} G(q_i)/\mathcal{U} \cong G(\Pi_{i\in \omega} {\mathbb F}_{q_i}/\mathcal{U}).$$
\end{theorem}

For $\Rightarrow$, Wilson first reduces to the case $G\equiv \Pi_{i\in \omega} S_i/\mathcal{U}$ (a non-principal ultraproduct of finite simple groups $S_i$). This uses a very nice observation of Felgner, that there is an $L_{\gp}$-sentence $\sigma$ which holds of every non-abelian simple groups, and with the property that any finite group $G$ satisfying $\sigma$ has non-abelian simple {\em socle} (the group generated by the minimal normal subgroups of $G$). The sentence $\sigma$ has form
$$\forall x \forall y\big[(x\neq 1 \wedge C_G(x,y)\neq 1) \to \bigcap_{g\in G}( C_G(x,y)C_G(C_G(x,y)))^g=1].$$

Wilson then analyses the possibilities for the $S_i$. It is easily seen that
$H=\Pi_{n\geq 5} \Alt(n)/\mathcal{U}$ is not simple, since finite alternating groups contain 3-cycles, and elements of increasingly large support, when written as products of 3-cycles, require increasingly many 3-cycles. The problem is that, naively,
$H$ might have an elementary substructure which is a simple group. To eliminate such possibilities, it suffices to show that, uniformly in $n$ , $\Alt(n)$ has an $\emptyset$-definable conjugacy-invariant family of elements of small support, and also such a family of increasingly large support, and elements of the latter cannot be written as a uniformly bounded product of elements of the former. Similar arguments work for groups elementarily equivalent to an ultraproduct of finite simple groups of increasingly large Lie rank -- that is, ultraproducts of groups $G_i$ of Lie type such that for each $n$, for almost all $i$ the group $G_i$ has Lie rank at least $n$.

The groups of Lie type each correspond to a Dynkin diagram. For twisted groups, such as ${}^2E_6(q)$, ${}^2F_4(q)$, etc., the Dynkin diagram has a symmetry which yields a `graph automorphism' of the corresponding untwisted group, essentially by permuting the root groups. One takes a product $\sigma$ of a graph automorphism and an appropriate `field automorphism' (arising from a power of the Frobenius), and, roughly speaking, takes the fixed points of $\sigma$ in the untwisted group (this description is not accurate -- see \cite[Chapter 13]{carter} for details.)

\medskip

Next, we consider soluble groups.

\begin{theorem} [Wilson \cite{wilson2}] \label{wil2}
There is an $L_{\gp}$-sentence $\sigma$ such that if $G$ is a finite group, then $G\models \sigma$ if and only if $G$ is soluble.
\end{theorem}
The sentence $\sigma$ asserts that there is no non-identity element $g$ which is a product of 56 commutators $[x,y]$ where each of $x,y$ is a conjugate of $g$. 

In one direction, it is clear that any soluble group satisfies this sentence $\sigma$ -- in fact, it satisfies the corresponding sentence with 56 replaced by any natural number. For suppose that $G$ is a group satisfying the above sentence $\sigma$, witnessed by $g\in G$. Let $N=\langle g\rangle^G$, the smallest normal subgroup of $G$ containing $g$. Then $g$ lies in the derived subgroup $N'$ of $N$, and hence $N'=N$, so $N$ is perfect and so not soluble, and hence $G$ is not soluble. In the other direction, Wilson uses Thompson's classification in \cite{thompson} of the {\em minimal} finite simple groups, that is, the minimal finite groups which are not soluble.

Note that it is {\em not} true that a pseudofinite group is soluble if and only if it satisfies $\sigma$. For if $G$ is an ultraproduct of a family of finite groups of increasingly large derived length then these groups satisfy $\sigma$ so by {\L}os's theorem $G\models \sigma$, but clearly $G$ is not soluble. 

Recall that the (soluble) {\em radical} $ R(G)$ of a group $G$ is the subgroup generated by the soluble normal subgroups of $G$. Always  $R(G)\triangleleft G$, and if $G$ is finite then $R(G)$ is soluble.

\begin{theorem} [Wilson\cite{wilson2}] \label{wilsonsol2}
There is an $L_{\gp}$-formula $\psi(x)$ such that if $G$ is a finite group then $\psi(G)=R(G)$.
\end{theorem}

The following questions appear to be open.

\begin{question} \rm 
1. Is there an $L_{\gp}$-sentence $\tau$ such that a finite group is nilpotent if and only if it satisfies $\tau$?

2. Is there an $L_{\gp}$ formula $\chi(x)$ which uniformly in finite groups defines the Fitting subgroup (the largest nilpotent normal subgroup)?
\end{question}

\section{Stability and generalisations, simple pseudofinite groups}

\subsection{Simple pseudofinite groups}

It follows fairly rapidly from the constructions of  the groups of Lie type, as described for example in \cite{carter}, that the finite groups of Lie type are uniformly definable in the corresponding finite fields, or, in the cases of Suzuki and Ree groups, in the corresponding difference fields. For the Suzuki and Ree groups this is noted in \cite{hrushrav}. In fact, we have the following. We say here that a class $\mathcal{C}$ of finite structures is uniformly definable (interpretable) in a class $\mathcal{D}$ if there are cofinite $\mathcal{C}'\subseteq \mathcal{C}$ and $\mathcal{D}'\subseteq \mathcal{D}$ and a bijection $f:\mathcal{D}'\to \mathcal{C}'$ such that for each $M\in \mathcal{D}'$, $f(M)$ is definable (respectively, interpretable) uniformly in $M$, i.e. always using the same formulas, but possibly allowing parameters. There is a corresponding notion of {\em uniform parameter bi-interpretability} -- for details see \cite{ryten}.

\begin{theorem}[Ryten]\label{rytenres}
Any family of finite simple groups of any fixed Lie type (other than Suzuki and Ree groups) is uniformly bi-interpretable (over parameters) with the corresponding family of finite fields.

(ii) The Ree groups ${}^2F_4(2^{2k+1})$ and the Suzuki groups ${}^2B_2(2^{2k+1})$ are uniformly parameter bi-interpretable with the difference fields $({\mathbb F}_{2^{2k+1}},x\mapsto x^{2^k})$, and the Ree groups ${}^2G_2(3^{2k+1})$ are uniformly parameter bi-interpretable with $({\mathbb F}_{3^{2k+1}}, x\mapsto x^{3^k})$.
\end{theorem}

Care is needed with the twisted groups. For example the unitary group PSU${}_n(q)$, which lives naturally as  a subgroup of PSL${}_n(q^2)$, is bi-interpretable (uniformly in $q$) with the field ${\mathbb F}_q$. It is a consequence of the main theorem of \cite{cdm} that ${\mathbb F}_q$ is {\em not} uniformly interpretable in ${\mathbb F}_{q^2}$. It follows that the groups
PSU${}_n(q)$ are {\em not} uniformly interpretable in the fields ${\mathbb F}_{q^2}$.

Extending remarks in Example~\ref{exstable}, we have 

\begin{theorem} \label{psfsimple}
(1) (Easy consequence of \cite{cdm}.) Any pseudofinite field has supersimple rank 1 theory.

(2) (From \cite{h}, resting on earlier work of Chatzidakis, Hrushovski and Peterzil (see \cite{ch} and \cite{chp}) Let $p$ be a prime, and let $m,n\in \omega$ with $m\geq 1, n>1$, and $(m,n)=1$. Let $\mathcal{C}_{m,n,p}$ be the class of finite {\em difference fields} (fields equipped with an automorphism) of form $({\mathbb F}_{p^{kn+m}}, {\rm Frob}^k)$ where $k\in \omega$.
Then any non-principal ultraproduct of $\mathcal{C}_{m,n,p}$ has supersimple rank 1 theory. 
\end{theorem}

In view of Theorem~\ref{rytenres} (2), this has particular interest for us in the cases $(m,n,p)=(1,2,2)$ and $(m,n,p)=(1,2,3)$. It yields the following.

\begin{corollary} [Hrushovski]\label{hrushryt}
Any simple pseudofinite group has supersimple finite rank theory.
\end{corollary}

\begin{remark} \rm 
Theorem~\ref{rytenres} was recently used by Nies and Tent \cite{nt} to show that

(1) finite simple groups are log-compressible, i.e., if $G$ is a finite simple group, there is a first order sentence $\phi$ in the language $L_{\gp}$, with unique model $G$, such that $\phi$ has length $O({\rm log} |G|)$, and more generally

(2) for any finite group $G$ there is such a sentence $\phi$ of length $O(({\rm log}|G|)^3)$.
\end{remark}

\begin{remark} \rm
The model theory of any non-principal ultraproduct $\Pi_{n\in {\mathbb N}} \Alt_n/\mathcal{U}$ is undecidable. The same holds for any non-principal ultraproduct $\Pi_{n\in {\mathbb N}} {\rm PSL}_n(q)/\mathcal{U}$. It seems hard to find a good reference, but see e.g. Section 6.3 of \cite{bunina}, or (originally)  Ershov \cite{ershov}.
\end{remark}

\subsection{Generalised stability for pseudofinite groups}

We aim here to give structural results for pseudofinite groups with stable, or more generally simple, or NIP, or NTP${}_2$, theory.

\begin{theorem}[\cite{mt}] \label{mtNIP}
(1) Let $\mathcal{C}$ be a class of finite groups such that all ultraproducts of members of $\mathcal{C}$ are NIP. Then there is $d\in \omega$ such that $|G:R(G)|\leq d$ for each $G\in \mathcal{C}$.

(2) If $G$ is pseudofinite NIP group with a fixed finite bound on the lengths of centraliser chains  then $G$ has an $\emptyset$-definable soluble normal subgroup of finite index.

(3) Any pseudofinite group with stable theory has an $\emptyset$-definable soluble normal subgroup of finite index.
\end{theorem}

\begin{remark} \label{NIPsol} \rm
1. In (2), the conclusion is false without some assumption like that on centralisers. Indeed, 
we give an example ({\em cf.} Theorem~\ref{mtNIP}(2)) of an NIP pseudofinite group which is not soluble-by-finite. Let $G={\rm SL}_2({\mathbb Z}_p)$, and for each $k>0$ let $G_k$ be the open normal subgroup of $G$ 
of form $$G_k:=\{\begin{pmatrix}1+a & b \\ c & 1+d \end{pmatrix}: a,b,c,d\in p^k{\mathbb Z}_p\},$$
a {\em congruence subgroup} of $G$. Then the groups $G_k$ are uniformly definable in the NIP valued field $({\mathbb Q}_p,+,\times)$ (since $k$ corresponds to an element of the value group, which is interpretable), so the quotients $G/G_k$ are uniformly interpretable finite groups. Let $\mathcal{U}$ be a non-principal ultrafilter on $\omega$, and put
$$G^*:=\Pi_{k\in \omega} (G/G_k)/\mathcal{U}.$$
Then $G^*$ is an NIP pseudofinite group. By $\omega_1$-saturation of ultraproducts, it has a normal subgroup $N$
such that $G^*/N\cong G$. In particular, $G^*$ is not soluble-by-finite.

2. Part (3) provides another route to the observation in the introduction that free groups are not pseudofinite. For by Sela's work they are known to be stable, and in the free non-abelian case they are clearly not soluble-by-finite.

\end{remark}

The proof of Theorem~\ref{mtNIP} makes essential use of Theorem~\ref{wilsonsol2}.

{\em Sketch Proof of Theorem~\ref{mtNIP}.} (1) Let $G\in \mathcal{C}$. Let $\psi(x)$ be as in Theorem~\ref{wilsonsol2}. For $G\in \mathcal{C}$ let $\bar{G}=G/R(G)$, and put $S:=\Soc(\bar{G})$
 (the direct product of the minimal normal subgroups). Then $S=T_1\times \ldots \times T_k$, where the $T_i$ are non-abelian finite simple groups.

{\em Claim 1.} There is a bound on $k$ as $G$ ranges through $\mathcal{C}$. Indeed, for each $i$ pick $x_i\in T_i\setminus Z(T_i)$ and $y_i\in T_i$ with $[x_i,y_i]\neq 1$. For $w\subset \{1,\ldots,k\}$ put $z_w=\Pi_{j\not\in w} y_j$. Then
$[x_j,z_w]=1 \Leftrightarrow j\in w$. Hence, the NIP assumption forces a bound on $k$.

{\em Claim 2.} There is a bound on the Lie rank of any $T_i$ (or on $t$ if $T_i=\Alt_t$). This is proved essentially as in Claim 1, as  otherwise some $T_i$ contains increasingly large direct powers of $\PSL_2$ or of $\Alt_4$.

{\em Claim 3.} The $T_i$ have bounded size. If this was false, then groups $G\in\mathcal{C}$ would contain arbitrarily large finite simple groups of fixed Lie rank (by Claims 1 and 2 and the classification of finite simple groups) so some ultraproduct would be a simple pseudofinite group, and (e.g. by Theorem~\ref{rytenres}) would interpret a pseudofinite field. But as noted in Examples~\ref{exstable}, pseudofinite fields do not have NIP theory.

By Claim 3, $|S|$ is bounded, and it follows easily that $|G:R(G)|$ is bounded.

(2) We may suppose  that $G=\Pi G_i/\mathcal{U}$ (an ultraproduct of finite groups), where each non-principal ultraproduct of the $G_i$ is elementarily equivalent to $G$. Thus by (1) there is a finite bound on $|G_i:R(G_i)|$. By stability of $G$  and the remarks before Theorem~\ref{psfsimple}, there is some $e\in\omega$ such that every centraliser chain in  $G$ has length at most $e$, and hence the same holds for any $G_i$. By a result of Kukhro \cite{kuk}, there is a function $f$ such that each group $R(G_i)$  has derived length at most $f(e)$. It follows that 
$R(G)$ (also defined by $\psi(x)$) is soluble.

\begin{example} \cite[Section 5]{mtstab} \rm 

(1) There is an $\omega$-stable pseudofinite group $G$ which is not nilpotent-by-finite. It has form $({\mathbb C},+)\rtimes \Gamma$ for some infinite but `small' $\Gamma \leq ({\mathbb C}^*,\cdot)$. This is a conglomeration of work of Chapuis, Simonetta, Khelif, and Zilber. The group has infinite Morley rank -- for Khelif has shown that any pseudofinite group of finite Morley rank is abelian-by-finite.

(2) The `Mekler construction' gives, for any odd prime $p$, examples of pseudofinite $\omega$-stable groups which are nilpotent of class 2 and exponent $p$ but not finite-by-abelian-by-finite. For the Mekler construction, see  \cite{hodges} or \cite{mekler}. The idea is to code graphs into nilpotent class 2 groups. Fix an odd prime $p$ and given a graph $\Gamma$ with vertex set $V$, let $G(\Gamma)$ be the group which is free nilpotent (on the generating set $V$) subject to being of class 2 and exponent 2, and  subject to the relations $[u,v]=1$ whenever  vertices $u,v$ are adjacent in $\Gamma$. Under reasonable conditions on $\Gamma$ (that it is a `nice graph'), properties such as stable and  simple are transferred from $\Gamma$ to $G(\Gamma)$ even though
$G(\Gamma)$ is not interpretable in $\Gamma$ (for simplicity, see \cite[Theorem 5.1]{baud}). Chernikov (personal communication) has shown that if $\Gamma$ is NIP then $G(\Gamma)$ is NIP, and it would be interesting to investigate which other model-theoretic conditions are preserved by the construction.

\end{example}

Next, we discuss pseudofinite groups with {\em simple} theory. Here, note that the examples (even supersimple of finite rank) include simple groups of Lie type over pseudofinite fields (by Corollary~\ref{hrushryt})
and also, for odd primes $p$, infinite extraspecial $p$-groups of exponent $p$, that is, groups $G$ of exponent  $p$ such that $G'=Z(G)=\Phi(G)\cong C_p$, where $\Phi(G)$ is the Frattini subgroup of $G$. Extraspecial $p$-groups have SU rank 1, and are finite-by-abelian but not abelian-by-finite. They have infinite descending chains of centralisers, and do not have a smallest finite index definable subgroup. For more detail see the Appendix of Milliet \cite{milliet}, or \cite[Proposition 3.11]{ms}.

The following result shows that ultraproducts of finite extraspecial groups are not simple unles at least one of the prime and the rank  (of the elementary abelian group $G/Z(G)$) is bounded.

\begin{proposition} 
For each $p,n\in {\mathbb N}^{>0}$ with $p$ prime, let $G_{p,n}$ be the extraspecial $p$-group of order $p^{2n+1}$, and let $\mathcal{U}$be an ultrafilter on the set of pairs $(p,n)$ such that for each $d\in {\mathbb N}$ there is $U\in \mathcal{U}$ such that for all $(p,n)\in U$ we have $p>d$ and $n>d$. Let $G:=\Pi_{p,n} G_{p,n}/\mathcal{U}$. 
Then $\Th(G)$ is not simple.
\end{proposition}

{\em Sketch Proof.} Let $Z:=Z(G)$ and $V=G/Z(G)$, an infinite abelian group. The commutator map  defines a non-degenerate bilinear map $\beta:V\times V\to Z$ given by $\beta(uZ,vZ)=u^{-1}v^{-1}uv$. Since this is definable in  $G$, and since $Z$ is the additive reduct of a pseudofinite field $K$, it follows from Granger \cite[Proposition 7.2.2]{granger} that $K$ is interpretable in $G$, as is the infinite-dimensional vector space structure of $V$ over $K$, and we may view $\beta$ as a definable symplectic bilinear form on $V$. By \cite[Proposition 7.4.1]{granger},  such structures are not simple.

\bigskip

Consider a class $\mathcal{C}$ of finite groups with all ultraproducts of $\mathcal{C}$ having simple theory. For $G\in \mathcal{C}$, $R(G)$ is uniformly $\emptyset$-definable (by \cite{wilson3}) and $\Soc(G/R(G))$ is a product of boundedly many non-abelian finite simple groups of bounded Lie rank, by variants of the proofs of Claims 1 and 2 above; see also Theorem~\ref{NTP2gp} below.

\begin{question} \rm In this setting, must $R(G)$ have bounded derived length, as $G$ ranges through $\mathcal{C}$?
\end{question}

If we assume that all ultraproducts of $\mathcal{C}$ are {\em supersimple}, then the answer is positive, by the following result of Milliet, a significant strengthening of results in \cite{ejmr}.

\begin{theorem} If $G$ is a pseudofinite group with supersimple theory, then $R(G)$ is definable and soluble
(and likewise, if we assume $G$ has finite SU-rank, then $\Fitt(G)$ is definable and nilpotent).
\end{theorem}

Thus, if $G$ is pseudofinite with superstable theory then $G$ has soluble radical $R(G)$, and if $S=\Soc(G/R(G))$, then $S=T_1\times \ldots \times T_k$ where the $T_i$ are non-abelian finite or pseudofinite simple groups. If $\bar{S}$ denotes the preimage of $S$ in $G$ then $G/\bar{S}$ embeds in $\Aut(T_1\times \ldots\times T_k)$.

We have not discussed  properties of SU-rank, but note that finite groups have SU-rank 0, and that if $G$ has supersimple theory of finite SU-rank and $H\leq G$ is definable, than $\SU(G)=\SU(H)+\SU(G/H)$, where $G/H$ denotes the interpretable set of left cosets of $H$ in $G$. 

The next result gives some information on supersimple pseudofinite groups of small SU-rank. Note the currently  essential use f the classification of finite simple groups (CFSG) in (3) -- it would be interesting to remove this. 

\begin{theorem} \label{ejmrsol}
Let $G$ be a pseudofinite group with supersimple theory, and assume that $T^{\eq}$ eliminates the quantifier $\exists^{\infty}$, where $T={\rm Th}(G)$. 

(1) \cite{er} If $SU(G=1$ then $G$ has a finite index definable characteristic subgroup $N$ such that $N'$ is a finite subgroup of $Z(N)$ (so $G$ is (finite-by-abelian)-by-finite). 

(2) \cite{ejmr} If $SU(G)=2$ then $G$ is soluble-by-finite.

(3) \cite{ejmr}  (CFSG)  If $G$ is a simple group and $SU(G)=3$  then $G \cong {\rm PSL}_2(K)$ for some pseudofinite field $K$.
\end{theorem}

Certain infinite (monomial) SU-rank versions of these results have recently been proved by Wagner. Parts (1) and (2) above are proved without the classification of finite simple groups. It should also be possible to remove the assumption on the quantifier $\exists^{\infty}$ here 
and also in Theorem~\ref{rankone} below. This assumption was natural in the context of \cite{ejmr} where the central context was that of groups with {\em measurable} theory, for which the assumption holds.

\medskip

We turn now to the NTP${}_2$ condition in the context of finite and pseudofinite groups. Here, we use the following consequence of the NTP${}_2$ condition.

\begin{lemma}\cite[Lemma 4.3]{mt3} \label{l:NTP2}
Let $G$ be an $\emptyset$-definable  group in a structure with  NTP${}_2$ theory, and $\psi(x,\bar{y})$ a formula implying $x\in G$. Then there is  $k=k_\psi\in {\mathbb N}$ such that the following holds. 
Suppose that $H$ is a subgroup of $G$, $\pi: H\longrightarrow \Pi_{i\in J} T_i$ is an epimorphism to the Cartesian product of the groups $T_i$, and
$\pi_j: H\longrightarrow T_j$ is for each $j\in J$ the composition of $\pi$ with the canonical projection $\Pi_{i\in J} T_i\to T_j$.
Suppose also that for each $j\in J$,  there is a subgroup $\bar{R}_j\leq G$  and group  $R_j <T_j$  with
$\bar{R}_j\cap H=\pi_j^{-1}(R_j)$, such that finite intersections of the groups $\bar{R}_j$ are uniformly definable by instances of $\psi(x,\bar{y})$.
Then $|J|\leq k$.
\end{lemma}

In the theorem below and its proof, we view $\Alt_n$ as having Lie rank $n$.

\begin{theorem}\label{NTP2gp}
Let $\mathcal{C}$ be a class of finite groups all of whose ultraproducts are NTP${}_2$. Then there is $d\in {\mathbb N}$ such that the following hold, where
 $G\in \mathcal{C}$ and $R(G)$ is  the soluble radical of $G$, with $\pi:G \to G/R(G)$ the natural map, and $S:=\Soc(G/R(G))$: the
 group $S$ is a direct product $T_1\times \ldots \times T_r$ of at most $d$ non-abelian simple groups $T_i$ which are of order at most 
$d$ or of Lie rank (possibly twisted) at most $d$, and $R(G)$ and the groups $\pi^{-1}(T_i)$ are uniformly definable, using
 finitely many formulas $\phi(x,\bar{y})$ as $G$ ranges through $\mathcal{C}$.
\end{theorem}

{\em Proof.} Using Wilson's Theorem~\ref{wilsonsol2}, we may suppose that $R(G)=1$ for $G\in \mathcal{C}$. 

{\em Claim 1.} $S$ is a direct product of a bounded number of simple groups.

{\em Proof of Claim 1.} Suppose that for each $e\in {\mathbb N}$ there is $G\in \mathcal{C}$ such that $S=T_1 \times \ldots \times T_m$ for $m\geq e$, where the $T_i$ are non-abelian simple groups. By \cite[Corollary 1.5]{ls} (together with the Feit-Thompson Theorem) there is a constant $c$ that that if $G$ is a finite non-abelian simple group  then every element of $G$ is a product of exactly $c$ conjugate involutions. In particular, there is $g=(g_1,\ldots,g_m)\in S$, where each $g_i$ has order 2, such that every element of $S$ is a product of $c$ conjugates of $g$. Since $S\triangleleft G$ it follows that $S$ is uniformly definable in $G$. 

Now, with $g$ as above and $I\subseteq \{1,\ldots,m\}$ let $g^{(I)}_j=g_j$ if $j\in I$, and $g^{(I)}_j=1$ otherwise. Put
$g^{(I)}:=(g^{(I)}_1,\ldots,g^{(I)}_m)$. Let $\pi_I$ be the projection of $S$ onto the coordinates indexed by $\{1,\ldots,m\}\setminus I$, and $S_I:=\Ker(\pi_I)$. Then the elements of $S_I$ are exactly the products  of at most $c$ conjugates in $S$ of $g^{(I)}$, so as $S$ is definable (uniformly as $G$ varies) so are the $S^{(I)}$. We may now apply the finitary version of Lemma~\ref{l:NTP2}, putting $H=S$ and 
$\bar{R}_i:=\pi_{\{i\}}^{-1}$ for each $i=\{1,\ldots,m\}$, to  conclude that some ultraproduct of $\mathcal{C}$ has TP${}_2$ theory.

\medskip

Given Claim 1, write $S=T_1 \times\ldots \times T_r$, where $r\leq d$ and the $T_i$ are non-abelian simple. It remains to prove

\medskip

{\em Claim 2.} There is a bound on the Lie rank of the $T_i$.

{\em Proof of Claim 2.} Since  the $T_i$ are uniformly definable, it suffices to show that any infinite ultraproduct of finite simple groups of increasingly large Lie rank has TP${}_2$ theory. We 
give a proof for alternating groups -- the proof for classical groups of Lie type is very similar and is only sketched here. So let $\mathcal{U}$ be an ultrafilter on ${\mathbb N}$ and
 $H:=\Pi_{n\in {\mathbb N}} \Alt_n/\mathcal{U}$. We view $\Alt_n$ as acting on $[n]:=\{1,\ldots,n\}$. It is well-known that the permutation group $(\Alt_n,[n])$ is uniformly definable in
 the abstract group $\Alt_n$. Likewise,  any $J\subset [n]$, is uniformly (in $n,J$) parameter-definable in  $\Alt_n$ as a set of form $\Fix(g)$ for appropriate $g$. Hence, subgroups of 
$\Alt_n$ of form $(\Alt_n)_{(J)}:=\{g\in \Alt_n: g|_J={\rm id}|_J\}$ are uniformly definable.
Now for increasingly large $m$ and $n>>m$, pick disjoint subsets $J_1,\ldots,J_m,K_{11},\ldots,K_{1m},\ldots,K_{m1},\ldots,K_{mm}$ of $[n]$ of size $m$. For each $i,j\in \{1,\ldots,m\}$ pick
 $a_{ij}\in \Alt_n$ with $J_i^{a_{ij}}=K_{ij}$. Also let $\psi(x,\bar{b}_i)$ be a formula defining $(\Alt_n)_{(J_i)}$. Let $\phi(x,\bar{b}_i a_{ij})$ be the formula expressing 
$x\in (\Alt_n)_{(J_i)}a_{ij}$. Then for each $i$ the formulas $\phi(x,\bar{b}_ia_{ij})$ are 2-inconsistent, and, essentially because of the disjointness of the $J_i$ and $K_{ij}$, for any
 $f:\{1,\ldots,m\}\to \{1,\ldots,m\}$, the set
$\phi(x,\bar{b}_i a_{i,f(i)}): 1\leq i\leq m\}$ is consistent. It follows by compactness that $\mathcal{C}$ has a T{P${}_2$ ultraproduct, a contradiction.

For the proof of Claim 2 when $H$ is an ultraproduct of classical groups there are several  arguments, and we omit some details. Suppose for example that $H$ is an ultraproduct of groups of the form 
PSL${}_{n_i}(q_i)$ where $n_i \to \infty$. By the argument in \cite[Proposition 3.11]{barbina}, there is a uniformly definable set $C$ of pairs $(g,g')$ of transvections in PSL${}_n(q)$ such that each pair determines a point of projective space, and a uniformly definable equivalence relation $E$ on $C$ such that $(g,g')E(h,h')$ if and only if $(g,g')$ and $(h,h')$ determine the same projective point. We may thus identify the corresponding projective space with $C/E$, with $G$ acting on it by conjugation. The argument then continues as in the last paragraph. For the symplectic, orthogonal, and unitary groups   similar results in \cite{barbina} can be applied.

\begin{remark} \rm
It follows from Theorem~\ref{NTP2gp} that if $G$ is a pseudofinite group with NTP${}_2$ theory then $G$ has an $\emptyset$-definable normal subgroup $R$ such that if $\bar{G}:=G/R$, then $\bar{G}$ has a definable normal subgroup $S$ (the group generated by the definable minimal normal subgroups of $\bar{G}$) which is a direct product of finitely many definable finite or pseudofinite simple groups. We do not know if $R$ must be soluble if $G$ has simple theory -- but note by Remark~\ref{NIPsol} (1) that in general $R$ need not be soluble, even assuming that $G$ has NIP theory.
\end{remark}

\subsection{Applications of generalised stability}

We discuss several ways in which the model theory of pseudofinite groups has potential applications in finite group theory, or at least provides a model-theoretic viewpoint. There is overlap with Section 4 of \cite{mt} and Section 6 of \cite{garcia}.

\medskip

{\bf 1. Indecomposability.} First, we mention a version of the well-known `Zilber Indecomposability Theorem' for groups of finite Morley rank, itself a generalisation of a classical result on algebraic groups. The result is due to Wagner  \cite[4.5.6]{wagner}, and the formulation below is in \cite[Remark 2.5]{er}.

\begin{theorem}[Indecomposability Theorem] \label{indec}
Let $G$ be a group interpretable in a supersimple finite $\SU$-rank theory, and let $\{X_i:i\in I\}$ be a collection of definable subsets of $G$. Then there exists a definable subgroup $H$ of $G$ such that:
 
 (i) $H\leq \langle X_i:i\in I\rangle$, and there are $n\in {\mathbb N}$, $\epsilon_1,\ldots,\epsilon_n\in \{-1,1\}$, and
 $i_1,\ldots,i_n \in I$, such that $H\leq X_{i_1}^{\epsilon_1}\ldots X_{i_n}^{\epsilon_n}$.
 
 (ii) $X_i/H$ is finite for each $i\in I$.
 
\noindent
If the collection of $X_i$ is setwise invariant under some group $\Sigma$ of definable automorphisms 
of $G$, then $H$ may be chosen to be $\Sigma$-invariant.
\end{theorem}

\begin{theorem}\cite[Theorem 4.2]{mt} \label{indecap}
Let $\mathcal{C}_\tau$ be the family of finite simple groups of fixed Lie type $\tau$ (possibly twisted), and let $\phi(x,\bar{y})$ be an $L_{\gp}$-formula. Then there is $d=d(\phi,\tau)$ such that
 if $G\in\mathcal{C}_\tau$, $\bar{a}\in G^{|\bar{y}|}$, and $X=\phi(G,\bar{a})$ satisfies $|X|>d$, then $G$ is a product of at most $d$ conjugates of $X\cup X^{-1}$.
\end{theorem}

There are analogues of Theorem~\ref{indec} already in \cite{hp}, for groups uniformly definable in finite fields -- see e.g. \cite[Proposition 1.13]{hp}. Various consequences are given there -- for example, in Proposition 4.3,  a new proof of a result of Nori on subgroups of GL${}_n(p)$ generated by elements of order $p$. A further application of such results is given by Lubotzky in \cite{lub}, in a proof of a result announced in \cite{kln}. Recall that, for $0<\epsilon\in {\mathbb R}$,   a finite $k$-regular graph $\Gamma$ with vertex set $V$ is called an {\em $\epsilon$-expander} if for every $A\subset V$ with $|A|\leq \frac{1}{2}|V|$ we have $\partial A|\geq \epsilon|A|$, where $\partial A$ is the set of vertices outside $A$ with a neighbour in $A$. Suzuki groups have also been shown to satisfy the theorem below -- see \cite{bgt}. 
\begin{theorem} \cite{kln}
There is $k\in {\mathbb N}$ and $0<\epsilon\in {\mathbb R}$ such that if $G$ is a finite simple group (not a Suzuki group), then  $G$ has a set of $k$ generators for which the Cayley graph Cay$(G,S)$ is an $\epsilon$-expander.
\end{theorem}

In the approach to this theorem in \cite{lub}, a key step is the following result. As explained in \cite{lub}, it follows almost immediately from \cite{hp}, or from Theorem~\ref{indec}. 

\begin{theorem} \cite[Theorem 4.1]{lub}
There is a function $f:{\mathbb N}\to {\mathbb N}$ such that if $G$ is a finite simple group of Lie type of rank $r$, but not of Suzuki type, then $G$ is a product of $f(r)$ copies of ${\rm SL}_2$.
\end{theorem}

\bigskip 

{\bf 2. Asymptotic classes.} 

The following notion was introduced by Elwes in \cite{elwes}, extending a 1-dimensional version explored in \cite{ms}.

\begin{definition} \label{asymp}\rm 
Let $\mathcal{C}$ be a class of finite $L$-structures. Then $\mathcal{C}$ is an {\em $N$-dimensional asymptotic class} if the following hold.

(i) For
every $L$-formula $\phi(\bar{x}, \bar{y})$ where $l(\bar{x})=n$ and $l(\bar{y})=m$,
there is a finite set of pairs $D \subseteq (\{0,\ldots,Nn\}
\times {\mathbb R}^{>0}) \cup \{(0,0)\}$ and for each $(d, \mu) \in
D$ a collection $\Phi_{(d,\mu)}$ of pairs of the form $(M,
\bar{a})$ where $M \in \mathcal{C}$ and $\bar{a} \in M^m$, so that $\{
\Phi_{(d, \mu)} : (d, \mu) \in D \}$ is a partition of 
$\{ (M, \bar{a}): M \in \mathcal{C},  \bar{a} \in M^m \}$, and
$$\big||\phi(M^n, \bar{a})| - \mu|M|^{\frac{d}{N}}\big|
=o(|M|^{\frac{d}{N}})$$ as $|M| \rightarrow \infty$ and $(M,
\bar{a}) \in \Phi_{(d,\mu)}$.

(ii) Each $\Phi_{(d, \mu)}$ is $\emptyset$-definable, that is to
say $\{ \bar{a} \in M^m : (M, \bar{a}) \in \Phi_{(d, \mu)}\}$ is
uniformly $\emptyset$-definable across $\mathcal{C}$.

\end{definition}

The class of all finite fields is, by the main theorem of \cite{cdm}, a {\em 1-dimensional asymptotic class} in the sense of \cite{ms}. Likewise, by \cite{ryten} the classes $\mathcal{C}_{1,2,2}$ and $\mathcal{C}_{1,2,3}$ of difference fields of form $({\mathbb F}_{2^{2k+1}},x\mapsto x^{2^k})$ and $({\mathbb F}_{3^{2k+1}},x\mapsto x^{3^k})$ respectively are 1-dimensional asymptotic classes. Elwes showed that if $\mathcal{C}$ and $\mathcal{C'}$ are families of finite structures and $f:\mathcal{C}\to \mathcal{C}'$ is a bijection such that for each $M\in \mathcal{C}$, $M$ and $f(M)$ are uniformly parameter-free bi-interpretable, then $\mathcal{C}$ is an asymptotic class if and only if $\mathcal{C}'$ is. With some additional work (because of use of parameters to interpret the fields in the groups), 
this yields

\begin{theorem}\cite{ryten} \label{rytenref}
Let $\mathcal{C}_\tau$ be the class of all finite simple groups of fixed Lie type $\tau$. Then $\mathcal{C}_\tau$ is an $N$-dimensional asymptotic class for some $N$ (and the values of $\mu$ in the definition are rational).
\end{theorem}

It is shown in \cite{ms} that if $M$ is an ultraproduct of an $N$-dimensional asymptotic class then $\Th(M)$ is supersimple of rank at most $N$. Furthermore, it is possible, using the definability clause (ii) in Definition~\ref{asymp},  consistently to assign a pair $(d,\mu)$ to every definable set so that certain basic counting axioms are satisfied in $M$; we say that $\Th(M)$ is {\em measurable}. It follows from Theorem~\ref{rytenref} and \ref{simplegroup} that any simple pseudofinite group has measurable theory in this sense. Measurability for groups is discussed further in \cite{ms}, \cite{er} and \cite{ejmr}, but not explored here.

\bigskip

We know that classes $\mathcal{C}_\tau$ of finite simple groups of fixed Lie type are uniformly definable in finite (difference) fields. In fact, much more is definable.  The asymptotic information in Theorem~\ref{rytenref} should have applications through the following result (see also Theorem~\ref{lmtmax} below). For the notion of restricted weight, see the discussion above \cite[Proposition 4.12]{lmt}. 

\begin{proposition} \cite[Proposition 4.12]{lmt} Let $\mathcal{C}_\tau$ be a class of finite simple groups $G(q)$ of fixed Lie type $\tau$, and let $V(\lambda)$ be an irreducible ${\mathbb F}_qG(q)$-module of restricted weight $\lambda$, with the action of $G(q)$ on $V(\lambda)$ given by $\rho(q)$. Then the structures $(G(q), V_\lambda(q),\rho(q))$ are uniformly definable in the fields ${\mathbb F}_q$ or in corresponding difference fields.
\end{proposition}

 The following result is proved in \cite[Proposition 2.2]{psw}, with the easy (ii) added in \cite[Theorem 4.7]{mt}. We do not give background on generic types for groups definable in simple theories, but refer to \cite{wagner} or \cite{psw}. In (ii), $G^{\circ}_M$ denotes the intersection of the $M$-definable subgroups of $G$ of finite index.

\begin{theorem}\label{3gens}
Let $T$ be a simple theory over a  countable language, $\bar{M}$ an $\omega_1$-saturated model of $T$ with a countable elementary substructure $M$, and $G$ an $\emptyset$-definable group in $\bar{M}$.
Let $p_1,p_2,p_3$ be three principal generic types of $G$ over  $M$. 

(i) There are $g_1,g_2\in \bar{M}$ such that $g_i\models p_i$ for $i=1,2$, $g_1$ and $g_2$ are forking-independent over $M$, and $g_1g_2\models p_3$.

(ii) If $r\in S_G(M)$ has realisations in $G_M^o$ then there are $a_i\in G$ with $a_i\models p_i$ (for $i=1,2,3$) such that $a_1a_2a_3\models r$. 
\end{theorem}

Using the asymptotic information in Theorem~\ref{rytenref}, this easily yields the following.

\begin{corollary}\label{pswcor}
Let $\mathcal{C}_\tau$ be as in Theorem~\ref{indecap}, and let $\phi_i(x,\bar{y})$ be formulas for $i=1,2,3$. Then there is $\mu\in {\mathbb Q}^{>0}$ such that for any sufficiently large $G\in \mathcal{C}_\tau$ and $\bar{a}_1,\bar{a}_2,\bar{a}_3\in G^{|\bar{y}|}$, if
$|\phi(G,\bar{a}_i)|\geq \mu |G|$ for each $i$, then
$$\phi_1(G,\bar{a}_1).\phi_2(G,\bar{a}_2).\phi_3(G,\bar{a}_3)=G.$$
\end{corollary}

The proof shows in addition that $|\frac{\phi_1(G,\bar{a}_1).\phi_2(G,\bar{a}_2)|}{|G|}\to 1$ as $|G|\to \infty$. We remark that the same result follows from Nikolov-Pyber \cite{np}, where it is rapidly derived from the following result of Gowers (and the Nikolov-Pyber result is about arbitrary sufficiently large subsets of $G$, not necessarily definable). 

\begin{proposition} (\cite{gowers}, see also \cite{np})
Let $G$ be a group of order $n$ such that the minimal degree of a nontrivial representation is $k$. If $A,B,C$ are three subsets of $G$ such that $|A|.|B|.|C|>\frac{n^3}{k}$, then there is $(a,b,c)\in A\times B\times C$ such that $ab=c$.
\end{proposition}

In particular, if $w(x_1,\ldots,x_d)$ is a non-trivial group word, then $w$ defines a map $G^d\to G$ by evaluation, and we denote the image of $w$ by $w(G)$. For example, if $w(x_1,x_2)=x_1^{-1}x_2^{-1}x_1x_2$ then $w(G)$ is the set of commutators of $GH$.  Using a result of Larsen \cite{larsen} (with an earlier version due to Borel) 
which says that in simple algebraic groups the word map is dominant, Corollary~\ref{pswcor} yields 

\begin{theorem} Let $w_1,w_2,w_3$ be non-trivial group words, and $\mathcal{C}_\tau$ a family of finite simple groups of fixed Lie type.
Then $w_1(G)w_2(G)w_3(G)=G$ for sufficiently large $G\in \mathcal{C}_\tau$. 
\end{theorem}

\begin{remark} \rm
1. There has been considerable recent literature on word maps, with much stronger results proved. For example, by \cite{lst}, if $w_1,w_2$ are non-trivial words, and $G$ is {\em any} sufficiently large finite simple group, then $w_1(G)w_2(G)=G$. For finite quasisimple groups (groups $G$ such that $G=G'$ and $G/Z(G)$ is non-abelian simple) this does not hold in general, but for any three non-trivial words we have $w_1(G)w_2(G)w_3(G)=G$ if $G$ is sufficiently large relative the the $w_i$ -- see \cite{lst2}. The famous Ore Conjecture states that if $G$ is a non-abelian finite simple group then every element of $G$ is a commutator. This has now been proved --
see \cite{lo}.

2. If $w(x_1,\ldots,x_d)$ is a group word, Theorem~\ref{rytenref} can be applied, within a family 
$\mathcal{C}_\tau$ of finite simple groups, to the formula $\phi(\bar{x},y)$ of form $w(x_1,\ldots,x_d)=y$, to give uniformity on the  asymptotic sizes of the preimages of the word map $w:G^d \to G$ for $G\in \mathcal{C}_\tau$.
\end{remark}

\bigskip

{\bf 3. Towards CFSG?} Given that pseudofinite simple groups have supersimple finite rank theory, one might (ambitiously) hope to classify them, under the additional assumption of supersimplicity, {\em without} using the classification of finite simple groups (CFSG). More generally, one might hope, without CFSG, to describe infinite families of finite simple groups all of whose ultraproducts are supersimple of finite rank. This is in the spirit of the Cherlin-Zilber Algebraicity Conjecture, which asserts that any simple group of finite Morley rank is isomorphic to a simple  algebraic group over an algebraically closed field. 

Parts (1) and (2) of Theorem~\ref{ejmrsol} are in this spirit.
One route in this direction would be to classify, without CFSG,  families of finite simple groups (with supersimple ultraproducts) with a BN pair. In a  major piece of work, Tits and Weiss \cite{tw} classified `Moufang' generalised polygons.  Dello Stritto used this to show that each of the parametrised families of  finite Moufang generalised polygons is an asymptotic class. This and further work of dello Stritto yields a description of groups with supersimple theory of finite SU-rank which have a definable spherical Moufang BN pair of Tits rank at least 2. 

\section{Pseudofinite permutation groups}

There are the beginnings of a structure theory of pseudofinite permutation groups, and of the model theory of families of finite permutation groups, in part under additional model-theoretic hypotheses. Recall that a permutation group $G$ on a set $X$ (here written $(G,X)$) is {\em primitive} if there is no proper non-trivial $G$-invariant equivalence relation on $X$, and is {\em definably primitive} if there is no proper non-trivial {\em definable} $G$-invariant equivalence relation on $X$. In finite permutation group theory, and to a lesser extent infinite permutation group theory, primitive permutation groups act as building blocks for all permutation groups, and many questions are reduced to problems on  primitive permutation groups. Mimicking a result from \cite{mp} in the finite Morley rank case, Elwes and Ryten used Theorem~\ref{indec} to prove the following.

\begin{proposition}
\label{defprim}
Let $(G,X)$ be a definably primitive permutation group definable in a supersimple finite rank theory $T$ such that $T^{\eq}$ eliminates $\exists^{\infty}$, and suppose that for $x\in X$ the point stabiliser 
$G_x$ is infinite. Then $G$ is primitive on $X$.
\end{proposition}

In a fundamental result, Hrushovski \cite{hrush} described possible definable transitive group actions on a strongly minimal set in a stable theory. Our nearest analogue in the pseudofinite case is the following, with the classification of finite simple groups currently needed in the description of case (3).

\begin{theorem}\label{rankone}
Let $(X,G)$ be a definably  primitive  pseudofinite permutation group in  a supersimple finite rank theory which eliminates $\exists^{\infty}$, and suppose that $\rk(X)=1$.
Let $S=\Soc(G)$.
Then one of the following holds. 

(i) $\rk(G)=1$, and $S$ is divisible torsion-free abelian or elementary abelian, has finite index in $G$, and acts regularly on $X$.  

(ii) $\rk(G)=2$. Here $S$ is abelian so regular and identified with $X$.
There is an interpretable pseudofinite field $F$ with additive group $X$, and $G\leq 
{\rm AGL}_1(F)$ (a subgroup of finite index), in the natural action.

(iii) $\rk(G)=3$. There is an interpretable pseudofinite field $F$, $S={\rm PSL}_2(F)$, 
${\rm PSL}_2(F)\leq G \leq {\rm P}\Gamma{\rm L}_2(F)$, and $X$ can be identified with  ${\rm PG}_1(F)$
in such a way that the action of $G$ on ${\rm PG}_1(F)$ is the natural one. 

\end{theorem}

We take the opportunity to fill a gap at the end of the proof of Theorem~\ref{rankone}, pointed out by Wagner. Right at the end of the proof of Lemma 5.15 of \cite{ejmr}, at the end of Section 5, it is asserted that $|G:{\rm PSL}_2(F)|$ is finite (to ensure $SU(G)=3$), and in particular if $B$ is a definable group of automorphisms of $F$ then $B$ is finite. The reason given is that otherwise there would be $b\in B$ such that $\Fix(b)$ is an infinite (definable) subfield of $F$, contradicting that $T$ has finite rank. This argument is not clear -- {\em a priori} all elements of $B$ could have the same finite fixed field. However, it can be shown that in such a case the orbits of $B$ on $F$ would be the classes of a definable equivalence relation on $F$ with infinitely many infinite classes, contradicting the assumption that $SU(F)=1$. 

In \cite{bc} Borovik and Cherlin answer a question first raised in \cite[Problem 1.6]{bt}, showing that there is a function $:{\mathbb N} \to {\mathbb N}$ such that if $(G,X)$ is a primitive permutation group of finite Morley rank then  $RM(G) \leq f(RM(X))$, where $RM$ denotes Morley rank. The proof uses the O'Nan-Scott-Aschbacher analysis of \cite{mp}, but, remarkably, though there is no classification of simple groups of finite Morley rank, uses many of the difficult tools developed with such a classification in mind. There is an analogous result for definably primitive permutation groups in o-minimal structures (where there is a classification of definably simple definable groups, due to Peterzil, Pillay, and Starchenko) in \cite{mmt}. In the pseudofinite case, we pose the following question. It is also raised in \cite{ejmr}, where Theorem 6.2 provides partial information, and material in \cite{lmt} should yield an answer.

\begin{problem} \rm Show that there is a function $f:{\mathbb N} \to {\mathbb N}$ such that if $(G,X)$ is a pseudofinite definably primitive definable permutation group in a supersimple theory of finite SU-rank then $SU(G) \leq f(SU(X))$.
\end{problem}

Given the rich literature on finite primitive permutation groups, it is natural to attempt to {\em classify} primitive  pseudofinite permutation groups. This was tackled in \cite{lmt}, with the main results from there sketched below. Recall first that a transitive permutation group $(G,X)$ is primitive if and only if each point stabiliser $G_x$ (for $x\in X$) is a maximal subgroup of $G$ -- in fact, the lattice of $G$-congruences on $X$ is naturally isomorphic to the lattice of groups between $G_x$ and $G$. If $G$ is transitive on $X$, then an {\em orbital graph} of $G$ on $X$ is a graph with vertex set $X$ and edge set some $G$-orbit on the set of unordered 2-element subsets of $X$. The following useful criterion for primitivity, due to D.G. Higman,  is well-known.

\begin{proposition} \cite{higman} Let $G$ be a transitive permutation group on a set $X$. Then the following are equivalent.

(1) $G$ is primitive on $X$,

(2) every orbital graph of $(G,X)$ is connected.
\end{proposition}

The following is now an elementary exercise.

\begin{proposition} \cite{lmt}
Let $(G,X)$ be an $\omega$-saturated transitive pseudofinite permutation group. Then the following are equivalent.

(1) $(G,X)$ is  primitive.

(2) If $x\in X$ then $G_x$ is boundedly maximal in $G$, that is, there is $d\in {\mathbb N}$ such that if $g,h\in G\setminus G_x$ then there are $x_1,\ldots,x_{d+1}\in G_x$ and $\epsilon_1,\ldots,\epsilon_d\in \{\pm 1\}$ such that
$h=x_1g^{\epsilon_1}x_2\ldots x_d g^{\epsilon_d}x_{d+1}$.

(3) There is $e\in {\mathbb N}$ such that each orbital graph of $(G,X)$ is connected of diameter at most $e$.
\end{proposition}

The following question was raised in \cite{lmt}.

\begin{question} \rm
Is there a primitive pseudofinite permutation group with infinite point stabiliser such that there is no finite bound on the diameters of the orbital graphs? By the last theorem, such a structure will not be $\omega$-saturated.
\end{question}

In \cite{lmt} a description, close to a full classification, is given of  primitive $\omega$-saturated pseudofinite permutation groups. It is involved, and we omit the details.

A key ingredient in \cite{lmt} is to
 consider pairs $(G,H)$ where $G$ is a finite simple group of Lie type and $H$ is a maximal (proper) subgroup of $G$ (named by a unary predicate). This is equivalent to considering the group $G$ together with a definable primitive action of $G$ on a set $X$, namely the set of left cosets of $H$ in $G$. (A permutation group $G$ on $X$ is {\em primitive} if there is no proper non-trivial $G$-invariant equivalence relation on $X$; this is equivalent to point stabilisers being maximal subgroups.) If $G=G(q)$ is a simple group of Lie type and $q=(q')^r$, then a {\em subfield subgroup} of $G$ is one of the form $G(q')$ (so of the same Lie type), embedded naturally. Such subgroups can be maximal if $r$ is prime.

\begin{theorem} \cite{lmt} \label{lmtmax}
Let $\tau$ be a fixed Lie type, and let $\mathcal{C}_{\tau,d}$ be the set of pairs $(G,H)$ where $G$ is a finite simple group of Lie type $\tau$, $H$ is a maximal subgroup of $G$, and if $H$ is a subfield subgroup then the corresponding field extension has degree at most $d$. Then 

(1) the class $\mathcal{C}_{\tau,d}$ is uniformly definable in the corresponding family of fields or difference fields, that is, there are finitely many tuples of formulas which serve (with suitable choice of parameters) to define all such pairs;

(2) any non-principal ultraproduct of such a family $\mathcal{C}_{\tau,d}$ will be a pair $(G^*,H^*)$ with supersimple finite rank theory, such that $H^*$ is maximal in $G^*$.
\end{theorem}

This theorem was mis-stated in \cite[Corollary 4.11]{lmt}, for the subgroups PSU$(n,q)$ are maximal but not uniformly definable in PSL$(n,q^2)$ -- see the comments after Theorem~\ref{rytenres} above. The pair $(G,H)$ is uniformly definable in the (difference) field, but not, in a few special cases such as this, in the larger group $G$.

The last assertion in (2) above  (maximality of $H^*$ in $G^*$)  follows from the remaining assertions, together with an argument using 
Theorem~\ref{indec}. This was used in \cite{lmt} to give a description of all $\omega$-saturated pseudofinite primitive permutation groups, that is $\omega$-saturated pseudofinite pairs $(G,H)$ with $H$ a maximal subgroup of $G$ which is {\em core-free}, that is, satisfies $\bigcap_{g\in G}H^g=\{1\}$.  Essentially, this is equivalent to describing families $\mathcal{F}_d$ of {\em finite} primitive permutation groups $G$ on sets $X$ such that, for every orbit $E$  of $G$ on the set $X^{[2]}$ of unordered 2-subsets of $X$, the graph on $X$ with edge set $E$ is connected of diameter at most $d$.

\section{Further directions}

The following well-known question, raised by Sabbagh,  has been open for a long time.

\begin{question} \rm Is there a finitely generated pseudofinite group?
\end{question}

For some discussion of this question, see Section 3 of \cite{ohp} (for example Proposition 3.9). The latter paper has  a number of interesting results on pseudofinite groups somewhat disjoint to this survey, such as the following analogue of the Tits Alternative.

\begin{theorem} \cite[Theorem 4.1]{ohp}
Let $G$ be an $\omega$-saturated pseudofinite group. Then either $G$ contains a free subsemigroup of rank 2, or $G$ is nilpotent-by-(uniformly locally finite).
\end{theorem}

Motivated by foundational questions in physics, Zilber \cite{zilber} has asked the following question. See also \cite{pillay-bohr} for a discussion of related topics on pseudofinite groups.

\begin{question}\rm (Zilber) Can an ultraproduct of finite groups have SO${}_3({\mathbb R})$ (or any compact simple real Lie group)  as a quotient? More generally, it would be interesting to identify positive sentences of $L_{\gp}$ which hold of all finite groups but not of all groups.
Here a sentence is positive if it is equivalent to one in prenex normal form with only the propositional connectives $\wedge$ and $\vee$; such sentences are preserved by group homomorphisms.
\end{question}

The following theorem answers a question raised in conversation by G. Levitt. The result may already be known.

\begin{theorem} \label{levittq}
(i) Let $S$ be the group of all permutations of a countably infinite set $X$. Then $S$  does not embed in any pseudofinite group.

(ii) There is a finitely generated group which does not embed in any pseudofinite group.
\end{theorem}

\medskip

{\em Proof.} (i) Let $\sigma$ be the sentence
$$\exists f\exists g \exists h([f^2,g]=1 \wedge [f,g]\neq 1\wedge h^{-1}fh=f^2).$$
Suppose $G\models \sigma$, with witnesses $f,g,h$. Then $C_G(f^2)>C_G(f)$ and $h^{-1}C_G(f)h=C_G(f^2)$. Hence $h$ has infinite order, so $G$ is infinite. Thus, if $H$ is a finite group, then $H\models \neg \sigma$, so every pseudofinite group satisfies $\neg \sigma$. As $\neg \sigma$ is universal, every group which embeds in a pseudofinite group satisfies $\neg \sigma$.

However, we claim that $S={\rm Sym}(X)\models \sigma$, where $X$ is a countably infinite set. Indeed, write $X=\bigcup_{i\in\omega} X_i$ as a disjoint union of infinite co-infinite subsets of $X$. For each $i$ put $X_i:=\{x_{ij}:j\in {\mathbb Z}\}$. Let $f$ act on $X$ by putting $f(x_{ij})= x_{i,j+1}$ for each $i,j$. Since $f$ and $f^2$ have the same cycle type (infinitely many infinite cycles and no other cycles) they are conjugate, that is, there is $h\in S$ with $h^{-1}fh=f^2$. Let $g$ be the element of $S$ such that $x_{0,2i}^g=x_{0,2i+2}$, with $g$ fixing all other elements of $X$.  Then $g\in C_G(f^2)\setminus C_G(f)$, as required.

(ii) The finitely generated subgroup $\langle f,g,h\rangle$ of $S$ also satisfies $\sigma$, so does not embed in any pseudofinite group.

\begin{remark} \rm If $G$ is a locally finite group then $G$ embeds in a pseudofinite group. Indeed, we may suppose $G$ is infinite. Let $\Delta$ be the atomic diagram of $G$, and $T$ be the theory of finite groups. Then clearly $T\cup \Delta$ is consistent, and any model of it is an infinite pseudofinite group which embeds $G$. This is a special case of a result of Malcev that  if $G$ is a group then $G$ embeds in some ultraproduct of  the finitely generated subgroups of $G$.
\end{remark}

Finally, in \cite{ba}, Bello Aguirre has begun an investigation into pseudofinite {\em rings}, by giving the following complete description of the generalised stability properties of pseudofinite residue rings. Similar results, but from a different viewpoint (quotients of prime ideals in non-standard elementary extensions of $({\mathbb Z},+,\times)$) have been obtained by D'Aquino and Macintyre. 

\begin{theorem} \cite{ba}
Let $\mathcal{U}$ be a non-principal ultrafilter on ${\mathbb N}$ and $F$ be the ring $\Pi_{n\in {\mathbb N}} {\mathbb Z}/n{\mathbb Z} /\mathcal{U}$. Then exactly one of the following holds, where $T=\Th(F)$.

(1) $T$ is NIP and there is a finite set $S$ of primes and some $U\in \mathcal{U}$ such that for $n\in U$, every prime divisor of $n$ lies in $S$.

(2) $T$ is supersimple of finite rank, and there is $d\in {\mathbb N}$ and $U\in \mathcal{U}$ such that each $n\in U$ is a product of at most $d$ prime powers, each with exponent at most $d$.

(3) $T$ is NTP${}_2$ but not simple or NIP, and there is $U\in \mathcal{U}$ and $d\in {\mathbb N}$ 
such that each $n\in U$ has at most $d$ prime divisors, but the conditions in (1) and (2) do not hold.

(4) $T$ is TP${}_2$, and for every $d\in {\mathbb N}$ there is $U=U_d\in \mathcal{U}$ such that each $n\in U$ has at least $d$ distinct prime divisors. 

\end{theorem}

The proof uses some model theory of $p$-adically closed fields and of their ultraproducts. A key point in (2) is that ${\mathbb Z}/p^d{\mathbb Z}$ is, for fixed $d$, uniformly (in $p$) coordinatised by ${\mathbb Z}/p{\mathbb Z}$. The proof of the TP${}_2$ condition in (4) uses Proposition~\ref{ntp2crit}. The arguments in case (2) have more recently been extended by Bello Aguirre to prove the following.

\begin{theorem} 
Let $d\in {\mathbb N}^{>0}$. Then the collection of all residue rings ${\mathbb Z}/p^d{\mathbb Z}$ forms a $d$-dimensional asymptotic class. 
\end{theorem}


\begin{thebibliography}{999}
\bibitem{ax} J. Ax, `The elementary theory of finite fields', Ann. Math. 88 (1968), 239--271.
\bibitem{bs} J. Baldwin, J. Saxl, `Logical stability in group theory', J. Austral. Math. Soc. 21 (1976), 267--276.
\bibitem{barbina} S. Barbina, `Reconstruction of classical geometries from their automorphism groups', J. London Math. Soc. (2) 75 (2007), 298--316.
\bibitem{baud} A. Baudisch, `Mekler's construction preserves CM-triviality', Ann. Pure Appl. Logic 115 (2002), 115--173.
\bibitem{ba} R. Bello Aguirre, `Generalised stability of ultraproducts of finite residue rings', arXiv:1503.04454.
\bibitem{bc} A. Borovik, G. Cherlin, `Permutation groups of finite Morley rank', in {\em Model theory with applications to algebra and analysis} (Eds. Z. Chatzidakis, H.D. Macpherson, A. Pillay, A.J. Wilkie), London Math. Soc. Lecture Notes 350, Cambridge University Press, Cambridge, pp. 59--124.
\bibitem{bt} A. Borovik. S. Thomas, `On generic normal subgroups', in {\em Automorphisms of first order structures} (Eds. R. Kaye, H.D. Macpherson), Clarendon Press, Oxford, 1994, pp. 319--324.
\bibitem{bgt} E. Breuillard, B. Green. T.C. Tao, `Suzuki groups as expanders', Groups Geom. Dyn. 5 (2011), 281--299.
\bibitem{bunina} E.I. Bunina, A.V. Mikhalev, `Elementary properties of linear groups and related problems', J. Math. Sci. 123 (2004),  3921--3985.
\bibitem{carter} R.W.  Carter, {\em Simple groups of Lie type}, Wiley, London, 1972.
\bibitem{casan} E. Casanovas, {\em Simple theories and hyperimaginaries}, Lecture Notes in Logic 39, Cambridge University Press, Cambridge, 2011.
\bibitem{cdm} Z. Chatzidakis, L. van den Dries, A.J. Macintyre, `Definable sets over finite fields', {\it J. Reine Angew. Math.} 
{\bf 427} (1992) 107--135.
\bibitem{ch} Z. Chatzidakis, E. Hrushovski, `The model theory of difference fields', 
Trans. Amer. Math. Soc. 351 (1999), 2997--3071.
\bibitem{chp} Z. Chatzidakis, E. Hrushovski, Y. Peterzil, `Model theory of difference fields II. Periodic ideals and the trichotomy in all characteristics', Proc. London Math. Soc. (3) 85 (2002), 257--311.
\bibitem{chern1} A. Chernikov, `Theories without the tree property of the second kind', Ann. Pure Appl. Logic, 165 (2014), 695--723.
\bibitem{cks} A. Chernikov, I. Kaplan, P. Simon, `Groups and fields with NTP2', Proc. Amer. Math. Soc. 143 (2015), 395--406.
\bibitem{elwes} R. Elwes, `Asymptotic classes of finite structures', J. Symb. Logic 72 (2007), 418--438.
\bibitem{ejmr} R. Elwes, E. Jaligot, H.D. Macpherson, M.J. Ryten, `Groups in supersimple and pseudofinite theories', Proc. London Math. Soc. (3) 103 (2011), 1049--1082.
\bibitem{er} R. Elwes, M. Ryten, `Measurable groups of low dimension', Math. Log. Quart. 54 (2008), 374--386.
\bibitem{ershov} J. Ershov, `Undecidability of the theories of symmetric and simple finite groups', Dokl. Akad. Nauk. SSSR 158 (1964), 777-779. 
\bibitem{garcia} D. Garcia, H.D. Macpherson, C. Steinhorn, `Pseudofinite structures and simplicity', J. Math. Logic 15 (2015) no. 1, 1550002.
\bibitem{gowers} T. Gowers, `Quasirandom groups', Combin. Prob. Comput. 17 (2008), 363--378.
\bibitem{granger} N. Granger, {\em Stability, simplicity and the model theory of bilinear forms}, Ph.D. Thesis, University of Manchester (1999). (www.maths.manchester.ac.uk/~mprest/)
\bibitem{hm} D. Haskell, H.D. Macpherson, `A version of o-minimality for the $p$-adics', J. Symb. Logic 62 (1997), 1075--1092.
\bibitem{higman} D.G. Higman, `Intersection matrices for finite permutation groups', J. Alg. 6 (1967), 22--42.
\bibitem{hodges} W. Hodges, {\em Model theory}, Cambridge University Press, Cambridge, 1993.
\bibitem{hrush} E. Hrushovski, `Almost orthogonal regular types', Ann. Pure Appl. Logic 45 (1989), 139--155.
\bibitem{hrushrav} E. Hrushovski, `Pseudofinite fields and related structures', in {\em Model theory and applications} (Eds. L. B\'elair, Z. Chatzidakis, P. d'Aquino, D. Marker, M. Otero, F. Point, A.J. Wilkie), Quaderni di Matematica vol. 11, Aracne, Rome, 2002, pp. 151--212.
\bibitem{h} E. Hrushovski, `The elementary theory of the Frobenius automorphisms', arXiv:math/0406514.
\bibitem{hp} A. Hrushovski, A. Pillay, `Definable subgroups of algebraic groups over finite fields', J. reine angew. Math. 462 (1995), 69--91.
\bibitem{kln} M. Kassabov, A. Lubotzky, N. Nikolov, `Finite simple groups as expanders', Proc. Natl. Acad. Sci. USA 103 (2006), no.16, 6116-6119.
\bibitem{kuk} E.I. Khukhro, `On solubility of groups with bounded centralizer chains', Glasgow Math. J. 51 (2009), 49--54.
\bibitem{kiefe} C. Kiefe, `Sets definable over finite fields: their zeta functions', Trans. Amer. Math. Soc. 223 (1976), 45--59.
\bibitem{kim} B. Kim, {\em Simplicity theory}, Oxford Logic Guides, Clarendon Press, Oxford, 2013.
\bibitem{larsen} M. Larsen, `Word maps have large image', Isr. J. Math. 139 (2004), 149--156.
\bibitem{lst} M. Larsen, A. Shalev, P. Tiep, `The Waring problem for finite simple groups', {\it Ann. Math.}  {\bf 174} (2011) 1885--1950.
\bibitem{lst2} M. Larsen, A. Shalev, P. Tiep, `The Waring problem for finite quasisimple groups', Int. Math. Res. Not. 2013 no.10,  2323--2348.

\bibitem{ls} M.W. Liebeck, A. Shalev, `Diameters of finite simple groups: sharp bounds and applications', Ann. Math. 154 (2001), 383--406.
\bibitem{lmt} M.W. Liebeck, H.D. Macpherson, K. Tent, `Primitive permutation groups of bounded orbital diameter', Proc. London Math. Soc. (3) 100 (2010), 216--248.
\bibitem{lo} M.W. Liebeck, E.A. O'Brien, A. Shalev, P.H. Tiep, `The Ore Conjecture', J. Euro. Math. Soc. 12 (2010), 939--1008.

\bibitem{lub} A. Lubotzky, `Finite simple groups of Lie type as expanders', J. Euro. Math. Soc. 13 (2011), 1331--1341.
\bibitem{mmt} H.D. Macpherson, A. Mosley, K. Tent, `Permutation groups in o-minimal structures', J. London Math. Soc. (2) 62 (2000), 650--670.
\bibitem{mp} H.D. Macpherson, A. Pillay, `Primitive permutation groups of finite Morley rank', Proc. London Math. Soc. (3) 70 (1995), 481--504.
\bibitem{ms} H.D. Macpherson, C. Steinhorn, `One-dimensional asymptotic classes of finite structures', {\it Trans. Amer. Math. Soc.} 
 360 (2008) 411--448.
\bibitem{mtstab} H.D. Macpherson, K. Tent, `Stable pseudofinite groups', J. Alg. 312 (2007), 550-561.
\bibitem{mt} H.D. Macpherson, K. Tent, `Pseudofinite groups with NIP theory and definability in finite simple groups', in {\em Groups and model theory} (Eds. L. Str\"ungmann, M. Droste, L. Fuchs, K. Tent), Contemp. Math. 576, Amer. Math. Soc., Providence, RI, 2012, 255-267.
\bibitem{mt3} H.D. Macpherson, K. Tent, `Profinite groups with NIP theory and $p$-adic analytic groups', arXiv:1603.02179.
\bibitem{mekler} A.H. Mekler, `Stability of nilpotent groups of class 2 and prime exponent', J. Symb. Logic 46 (1981), 781--788.
\bibitem{milliet} C. Milliet, `Definable envelopes in groups having a simple theory', https://hal.archives-ouvertes.fr/hal-00657716v2/document
\bibitem{np} N. Nikolov, L. Pyber, `Product decompositions of quasirandom groups and a Jordan type theorem', 
J. Euro. Math. Soc., 13 (2011) 1063--1077
\bibitem{nt} A. Nies, K. Tent, `Describing finite groups by short first-order sentences', arXiv:1409:8390
\bibitem{ohp} A. Ould Houcine, F. Point, `Alternatives for pseudofinite groups', J. Group Theory 16 (2013), 461--495.
\bibitem{pillay} A. Pillay, {\em Geometric stability theory}, Oxford Logic Guides, Clarendon Press, Oxford, 1996.
\bibitem{psw} A. Pillay, T. Scanlon, F. Wagner, `Supersimple fields and division rings', Math. Research Letters 5 (1998) 473-483. 
\bibitem{pillay-bohr} A. Pillay, `Remarks on compactifications of pseudofinite groups', arXiv:1509.02895.
\bibitem{point} F. Point, `Ultraproducts and Chevalley groups', Arch. Math. Logic 38 (1999),  355-372.
\bibitem{ryten} M.J. Ryten, {\em Model theory of finite difference fields and simple groups}, Ph.D. Thesis, University of Leeds (2007).\hfill 
\phantom{hi}
(http://www1.maths.leeds.ac.uk/Pure/staff/macpherson/ryten1.pdf)
\bibitem{shalev1} A. Shalev, `Characterisation of $p$-adic analytic groups in terms of wreath products', J. Algebra 145 (1992), 204--208.
\bibitem{simon} P. Simon, {\em A guide to NIP theories}, Lecture Notes in Logic 44, Cambridge University Press, Cambridge, 2015.
\bibitem{tz} K. Tent, M. Ziegler, {\em A course in model theory}, Lecture Notes in Logic 40, Cambridge University Press, Cambridge, 2012.
\bibitem{thompson} J.G. Thompson, `Nonsolvable finite groups all of whose local subgroups are solvable (Part I)', Bull. Amer. Math. Soc. (NS) 74 (1968), 383--437.
\bibitem{tw} J. Tits, R. Weiss, {\em Moufang polygons}, Springer monographs in mathematics, Springer, Berlin, 2002.
\bibitem{ugurlu} P. Ugurlu, `Pseudofinite groups as fixed points of simple groups of finite Morley rank', J. Pure Applied Alg. 217 (2013), 892--900.
\bibitem{wagner} F. Wagner, {\em Stable groups}, London Math. Soc. Lecture Notes no. 240, Cambridge University Press, 1997.
\bibitem{wagner2} F. Wagner, {\em Simple theories}, Kluwer, Dordrecht, 2000. 
\bibitem{wilson1} J.S. Wilson, `On pseudofinite simple groups', J. London
Math. Soc. (2) 51 (1995),
 471--490.
\bibitem{wilson2} J.S. Wilson, `Finite axiomatisation of finite soluble groups', J. London Math. Soc. 74 (2006), 566--582.
\bibitem{wilson3} J.S. Wilson, `First-order characterization of the radical of a finite group', J. Symb. Logic 74 (2009),  1429--1435.
\bibitem{zilber} B. Zilber, `Perfect infinities and finite approximation', {\em Infinity and truth} 199--223, Lect. Notes Ser. Inst. Math. Sci. Natl. Univ. Singap., 25, World Sci. Publ., Hackensack, NJ 2014.
\end{thebibliography}
\end{document}